\documentclass[12pt]{article}

\usepackage{fullpage}
\usepackage{latexsym}
\usepackage{amsmath}
\usepackage{amsfonts}
\usepackage{authblk}
\usepackage{xcolor}
\usepackage[numbers]{natbib}

\newtheorem{Definition}{Definition}[section]
\newtheorem{Theorem}{Theorem}

\newtheorem{Lemma}{Lemma}
\newtheorem{Corollary}{Corollary}

\newenvironment{Proof}{\noindent{\bf Proof: }}{\qed}

\newcommand{\qed}{\protect\raisebox{.3ex}{\framebox[.5em]{\rule{0em}{.6ex}}}}

\newcounter{listCounter}

\newcommand{\ListLengths}{\setlength{\itemsep}{0ex}\setlength{\topsep}{1ex}\setlength{\partopsep}{0ex}}

\title{Noncommutativity in the analysis of piecewise discrete-time dynamical systems}
\author[a]{Benjamin T. Hendel\thanks{Corresponding author}}
\author[b]{Rafael Ruggiero}
\affil[a]{Independent Researcher, bthendel@gmail.com, Fort Worth, Texas, USA} 
\affil[b]{Independent Researcher, ruggiero3n1@gmail.com, São Paulo, Brazil}

\begin{document}
\maketitle

\begin{abstract}
In this paper, we present a new method for the analysis of piecewise dynamical systems that are similar to the Collatz conjecture in regard to certain properties of the commutator of their sub-functions. We use the fact that the commutator of polynomials $E(n)=n/2$ and $O(n)=(3n+1)/2$ is constant to study rearrangements of compositions of $E(n)$ and $O(n)$. Our main result is that for any positive rational number $n$, if $(E^{e_1} \circ O^{o_1} \circ E^{e_2} \circ \dotsb \circ O^{o_l} \circ E^{e_{l+1}})(n)=1$, then $(E^{e_1} \circ O^{o_1} \circ E^{e_2} \circ \dotsb \circ O^{o_l} \circ E^{e_{l+1}})(n) = \lceil(E^{e_1 + \dotsb + e_{l+1}} \circ O^{o_1 + \dotsb + o_{l}})(n)\rceil$, where exponentiation is used to denote repeated composition and $e_i$ and $o_i$ are positive integers. Composition sequences of this form have significance in the context of the Collatz conjecture. The techniques used to derive this result can be used to produce similar results for a wide variety of repeatedly composed piecewise functions.
\end{abstract}


\section{Introduction}

The analysis of piecewise discrete-time systems is of considerable interest due to the fact that these systems can accurately model a wide variety of natural phenomena. Applications of piecewise discrete-time linear systems range from numerically simulating differential equations with boundary conditions (with the finite difference method, for instance) \citep{Iserles} to the study of real-world control systems \citep{Feng}. Naturally, the study of discrete-time systems leads to questions about these systems' long-term equilibria \citep{Lin}. It is painfully clear that Mathematics in its current state is ill-equipped to answer general questions about the limiting behavior of discrete dynamical systems. This is evidenced by the field's inability to prove or disprove the Collatz conjecture, which is based on a very simple observation concerning a deceptively simple toy model of a piecewise discrete-time dynamical system. Let $E(n)=n/2$ and $O(n)=(3n+1)/2$, where $n$ is a positive rational number. Define
\[   
f(n) = 
\begin{cases}
E(n) &\quad\text{if $n$ is even}\\
O(n) &\quad\text{if $n$ is odd}\\
\end{cases}
\]
where the domain of $f$ is restricted to the positive integers. The Collatz conjecture states that for any $n \ge 3$, there exists some $w$ such that $f^w(n)\equiv(f \circ \dotsb \circ f)(n)=1$ \citep{Mimuro}. The conjecture has been numerically verified to hold up to $n=2^{68}$ \citep{Barina}, but a proof remains elusive. In this paper, we will present a new method for the analysis of piecewise dynamical systems that are similar to the Collatz conjecture in regard to certain properties of the commutator of their sub-functions ($E$ and $O$ in the case of $f$).

It is always the case that either $f(n)=E(n)$ or $f(n)=O(n)$, so it is always possible to write $f^w(n)=(E^{e_1} \circ O^{o_1} \circ E^{e_2} \circ O^{o_2} \circ \dotsb \circ O^{o_j} \circ E^{e_{l+1}})(n)$ for some natural numbers $e_1, \dots, e_{l+1}$ and $o_1, \dots, o_l$. We will study rearrangements of $E$ and $O$ compositions and their relationship to repeated compositions of $f$. In particular, we find in Corollary 1 that
when $(E^{e_1} \circ O^{o_1} \circ E^{e_2} \circ O^{o_2} \circ \dotsb \circ O^{o_l} \circ E^{e_{l+1}})(n)=1$ and $e_1, \dots, e_{l+1}$ and $o_1, \dots, o_l$ are positive integers, it is the case that
\[
(E^{e_1} \circ O^{o_1} \circ E^{e_2} \circ \dotsb \circ O^{o_l} \circ E^{e_{l+1}})(n) = \lceil(E^{e_1 + \dotsb + e_{l+1}} \circ O^{o_1 + \dotsb + o_{l}})(n)\rceil.
\]
Composition sequences of this form are significant in the context of the Collatz conjecture because it is trivial to pair the orbits of most even numbers under $f$ with them.

Theorem 1 is a generalization of this result that holds when $(E^{e_1} \circ O^{o_1} \circ E^{e_2} \circ O^{o_2} \circ \dotsb \circ O^{o_l} \circ E^{e_{l+1}})(n)=W$ for some positive rational number W. This result is derived by using properties of the commutator $[O^a, E^b] \equiv O^a \circ E^b - E^b \circ O^a$. Specifically, we make use of the fact that the value of $[O^a, E^b]$ is constant for positive $a$ and $b$ values.  The techniques applied here can be used in a more general setting to analyze repeatedly composed piecewise functions with similarly behaved sub-functions. 

\section {Some properties of $O(n)$ and $E(n)$}
\begin{Lemma}
	For positive integers $w$, $E^w(n)=\left(\frac{1}{2}\right)^w n$ and \[O^w(n)=\left(\frac{3}{2}\right)^w (n+1)-1.\]
\end{Lemma}

\begin{Proof}
For repeated applications of $E(n)$ the result is trivial. For repeated applications of $O(n)$, the result follows from induction on $w$.
\end{Proof}

\vspace{.25 in}

\begin{Definition}
	We define the commutator of functions $A(x)$ and $B(x)$ as $[A,B](x) \equiv A \circ B - B \circ A$.
\end{Definition}

\begin{Lemma}
	For positive integers $a$ and $b$\, the commutator $[O^a, E^b](n)$ has a constant positive value for all positive rational numbers $n$. 
\end{Lemma}

\begin{Proof}
We can use Lemma 1 to directly compute
\begin{align*}
[O^a,E^b](n) &\equiv (O^a \circ E^b)(n)-(E^b \circ O^a)(n)\\
&=\left(\frac{3}{2}\right)^a \left(\frac{n}{2^b}+1\right)-1 - \frac{1}{2^b} \left(\frac{3}{2}\right)^a (n+1) + \frac{1}{2^b}\\
&=\left(\frac{3}{2}\right)^a\left(1-\frac{1}{2^b}\right)+\frac{1}{2^b}-1
\end{align*}
so $[O^a,E^b]$ is constant since $n$ vanishes. To show that $[O^a,E^b]>0$, note that
\begin{align*}
[O^a,E^b]>0 &\Leftrightarrow\left(\frac{3}{2}\right)^a\left(1-\frac{1}{2^b}\right)+\frac{1}{2^b}-1 > 0\\
&\Leftrightarrow\left(\frac{3}{2}\right)^a\left(1-\frac{1}{2^b}\right)+\frac{1}{2^b} > 1\\
&\Leftrightarrow\left(\frac{3}{2}\right)^a\left(2^b-1\right)+1 > 2^b.\\
\end{align*}
Since the smallest possible value of $a$ and $b$ is 1, we have
\[
\left(\frac{3}{2}\right)^a\left(2^b-1\right)+1 \ge \frac{3}{2}\left(2^b-1\right)+1>2^b
\]
and the result is proven.
\end{Proof}
\vspace{.25 in}

Since we have shown that the value of $[O^a, E^b](n)$ does not depend on $n$, we will adopt the notation $[O^a, E^b]$ to denote the commutator when applicable.

\begin{Lemma}
	Suppose that $G(n) = (E^{e_1} \circ O^{o_1} \circ \dotsb \circ E^{e_j} \circ O^{o_j})(n)$ for natural numbers $e_1, \dots, e_j$ and $o_1, \dots, o_j$. Then, for positive integer $w$ and real number $H$,
	\[O^w(G(n)+H) = (O^w \circ G)(n) + \left( \frac{3}{2}\right)^w H\]
	and
	\[E^w(G(n)+H) = (E^w \circ G)(n) + \left( \frac{1}{2}\right)^w H.\]
\end{Lemma}

\begin{Proof}
	Using Lemma 1 we have
	\begin{align*}
	O^w(G(n)+H) &= \left(\frac{3}{2}\right)^w (G(n)+H+1)-1 \\
	&= \left(\frac{3}{2}\right)^w G(n)+ \left(\frac{3}{2}\right)^w H + \left(\frac{3}{2}\right)^w -1 \\
    &= (O^w \circ G)(n) + \left( \frac{3}{2}\right)^w H
	\end{align*}
	and
	\[
	E^w(G(n)+H) = \frac{1}{2^w} (G(n)+H) = \frac{1}{2^w} G(n) + \frac{1}{2^w} H = (E^w \circ G)(n) + \left( \frac{1}{2}\right)^w H.
	\]
\end{Proof}

\begin{Lemma}
	Suppose that $n$ is a positive rational number and that $(E^{e_1} \circ O^{o_1} \circ E^{e_2} \circ \dotsb \circ O^{o_l} \circ E^{e_{l+1}})(n) = W$ for positive rational number $W$ and positive integers $e_1, \dots, e_{l+1}$ and $o_1, \dots, o_l$. Let $\sigma_e = e_1+ \dots+ e_{l+1}$ and $\sigma_o = o_1+ \dots+ o_{l}$. Then
	\[
	(E^{e_1} \circ O^{o_1} \circ E^{e_2} \circ \dotsb \circ O^{o_l} \circ E^{e_{l+1}})(n) = (E^{\sigma_e} \circ O^{\sigma_o})(n)+C=W.
	\]
	where
	\[C=W - \frac{3^{\sigma_o} (n+1) -2^{\sigma_o}}{2^{\sigma_o+\sigma_e}}.\]
	
\end{Lemma}
\begin{Proof}
	We can compute the value of $C$ directly by noting that
	\[
	C = W - (E^{\sigma_e} \circ O^{\sigma_o})(n)
	\]
	and applying Lemma 1 to get
	\[
    C =W - \frac{1}{2^{\sigma_e}} \left[ \left(\frac{3}{2}\right)^{\sigma_o} (n+1)-1 \right] =W - \frac{3^{\sigma_o} (n+1) -2^{\sigma_o}}{2^{\sigma_o+\sigma_e}}.
    \]
\end{Proof}
\vspace{.25 in}

\begin{Lemma}
Suppose that $n$ is a positive rational number and that $(E^{e_1} \circ O^{o_1} \circ \dotsb \circ E^{o_l} O^{o_l})(n) = W$ for positive rational number $W$ and positive integers $e_1, \dots, e_{l+1}$ and $o_1, \dots, o_l$. Let $\sigma_e = e_1 + \dotsb + e_{l}$. Then, for $l > 1$,
\[
(E^{e_1} \circ O^{o_1} \circ \dotsb \circ E^{e_{l}} \circ O^{o_l})(n) = (E^{e_1 + \dotsb + e_{l}} \circ O^{o_1 + \dotsb + o_l})(n)+C=W.
\]
where
\[C=\sum_{i=1}^{l-1}\left(\frac{1}{2}\right)^{\sigma_e-\zeta (i,l)}\left(\frac{3}{2}\right)^{\gamma(i,l)}[O^{o_i},E^{\zeta (i,l)}]\]
and $\zeta(i,l) = \sum_{l \ge j > i} e_j$ and $\gamma(i,l) = \sum_{j<i\leq l} o_j$.
\end{Lemma}
\begin{Proof}
Let $\sigma_o = o_1 + \dotsb + o_{l}$. For any choice of $e_1,\dots,e_{l+1}$ and $o_1, \dots, o_l$ where $l \ge 2$ we can use Lemma 3 and the definition of the commutator to get that
\begin{align*}
(E^{\sigma_e} \circ O^{\sigma_o})(n) &= (E^{e_1} \circ E^{\sigma_e - e_1} \circ O^{o_1} \circ O^{\sigma_o - o_1})(n)\\&=E^{e_1}((E^{\sigma_e - e_1} \circ O^{o_1})(O^{\sigma_o - o_1}(n)))\\&=E^{e_1}((O^{o_1} \circ E^{\sigma_e - e_1})(O^{\sigma_o - o_1}(n))-[O^{o_1},E^{\sigma_e - e_1}])\\&=(E^{e_1} \circ O^{o_1} \circ E^{\sigma_e - e_1} \circ O^{\sigma_o - o_1})(n)-\left(\frac{1}{2}\right)^{e_1}[O^{o_1},E^{\sigma_e - e_1}].
\end{align*}
The same procedure can be carried out for $(E^{\sigma_e - e_1} \circ O^{\sigma_o - o_1})(n)$ if $l \ge 3$ and the result can be substituted into the equation above. This procedure can be continued for $l-1$ total steps to produce $(E^{e_1} \circ O^{o_1} \circ \dotsb \circ E^{e_{l}} \circ O^{o_l})(n)$. If we do not apply Lemma 3 and only use the definition of the commutator it follows that
\begin{align*}
(E^{\sigma_e} \circ O^{\sigma_o})(n) =E^{e_1}((O^{o_1} \circ E^{\sigma_e - e_1}\circ O^{\sigma_o - o_1})(n)-[O^{o_1},E^{\sigma_e - e_1}])\\=E^{e_1}(O^{o_1}( E^{e_2} \circ O^{o_2} \circ E^{\sigma_e - e_1 - e_2} \circ O^{\sigma_o - o_1})(n)-[O^{o_2},E^{\sigma_e - e_1 - e_2}])-[O^{o_1},E^{\sigma_e - e_1}])\\=E^{e_1}(O^{o_1}(E^{e_2}(O^{o_2}(\dotsb E^{e_{l-1}}((O^{o_{l-1}} \circ E^{e_{l}})( O^{o_l}(n))\\-[O^{o_{l-1}},E^{e_{l}}])-[O^{o_{l-2}},E^{e_{l-1}+ e_{l-2}}])-\dotsb)-[O^{o_1},E^{\sigma_e-e_1}]).
\end{align*}
Applying Lemma 3 to the above expression from the inside out yields the result.
\end{Proof}
\vspace{.25 in}

\section {Results}

\begin{Theorem}
Suppose that $n$ is a positive rational number and that $(E^{e_1} \circ O^{o_1} \circ E^{e_2} \circ \dotsb \circ O^{o_l} \circ E^{e_{l+1}})(n) = W$ for positive rational number $W$ and positive integers $e_1, \dots, e_{l+1}$ and $o_1, \dots, o_l$. Let $\sigma_e = e_1+ \dots+ e_{l+1}$ and $\sigma_o = o_1+ \dots+ o_{l}$. Then, for $l>1$,
\[
(E^{e_1} \circ O^{o_1} \circ E^{e_2} \circ \dotsb \circ O^{o_l} \circ E^{e_{l+1}})(n) = (E^{\sigma_e} \circ O^{\sigma_o})(n)+C=W.
\]
where $0 < C < W$.
\end{Theorem}
\begin{Proof}
From Lemma 4 we have
\[
C=W - \frac{3^{\sigma_o} (n+1) -2^{\sigma_o}}{2^{\sigma_o+\sigma_e}}<W.
\]
Applying Lemma 5 we get that
\begin{align*}
(E^{e_1} \circ O^{o_1} \circ E^{e_2} \circ \dotsb \circ O^{o_l})(E^{e_{l+1}}(n))\\= (E^{e_1} \circ O^{o_1} \circ E^{e_2} \circ \dotsb \circ O^{o_l} \circ E^{e_{l+1}})(n)\\=(E^{e_1 + \dotsb + e_{l}} \circ O^{o_1 + \dotsb + o_l})(E^{e_{l+1}}(n))+\sum_{i=1}^{l-1}\left(\frac{1}{2}\right)^{\sigma_e -\zeta (i,l)}\left(\frac{3}{2}\right)^{\gamma(i,l)}[O^{o_i},E^{\zeta (i,l)}]\\=W
\end{align*}
for any positive integer $e_{l+1}$. Applying the definition of the commutator and Lemma 3 we get that
\begin{align*}
(E^{e_1 + \dotsb + e_{l}} \circ O^{o_1 + \dotsb + o_l})(E^{e_{l+1}}(n)) &= (E^{e_1 + \dotsb + e_{l}} \circ O^{o_1 + \dotsb + o_{l}}\circ E^{e_{l+1}})(n)\\&= E^{e_1 + \dotsb + e_{l}} ((O^{o_1 + \dotsb + o_{l}}\circ E^{e_{l+1}})(n))\\&=E^{e_1 + \dotsb + e_{l}} ((E^{e_{l+1}} \circ O^{o_1 + \dotsb + o_{l}})(n)+[O^{o_1 + \dotsb + o_l},E^{e_{l+1}}])\\&=(E^{e_1 + \dotsb + e_{l+1}} \circ O^{o_1 + \dotsb + o_{l}})(n)\\&+\left(\frac{1}{2}\right)^{e_1 + \dotsb + e_l} [O^{o_1 + \dotsb + o_l},E^{e_{l+1}}]
\end{align*}
so
\[C=\left(\frac{1}{2}\right)^{e_1 + \dotsb + e_l} [O^{o_1 + \dotsb + o_l},E^{e_{l+1}}]+\sum_{i=1}^{l-1}\left(\frac{1}{2}\right)^{\sigma_e -\zeta (i,l)}\left(\frac{3}{2}\right)^{\gamma(i,l)}[O^{o_i},E^{\zeta (i,l)}].\]
We have
\[
\left(\frac{1}{2}\right)^{e_1 + \dotsb + e_l} [O^{o_1 + \dotsb + o_l},E^{e_{l+1}}] > 0
\]
and for each of the terms of the sum
\[
\left(\frac{1}{2}\right)^{\sigma_e -\zeta (i,l)}\left(\frac{3}{2}\right)^{\gamma(i,l)}[O^{o_i},E^{\zeta (i,l)}] > 0
\]
since Lemma 2 implies that every commutator term is greater than zero. So $C > 0$.
\end{Proof}
\vspace{.25 in}

\begin{Corollary}
Suppose that $n$ is a positive even number and that $(E^{e_1} \circ O^{o_1} \circ E^{e_2} \circ \dotsb \circ O^{o_l} \circ E^{e_{l+1}})(n) = 1$ for positive integers $e_1, \dots, e_{l+1}$ and $o_1, \dots, o_l$. Let $\sigma_e = e_1+ \dots+ e_{l+1}$ and $\sigma_o = o_1+ \dots+ o_{l}$. Then, for $l>1$,
\[
(E^{e_1} \circ O^{o_1} \circ E^{e_2} \circ \dotsb \circ O^{o_l} \circ E^{e_{l+1}})(n) = \lceil(E^{\sigma_e} \circ O^{\sigma_o})(n)\rceil.
\]
\end{Corollary}
\begin{Proof}
The result follows from Theorem 1.
\end{Proof}
\vspace{.25 in}

\begin{Corollary}
Suppose that $n$ is a positive even integer and that $(E^{e_1} \circ O^{o_1} \circ E^{e_2} \circ \dotsb \circ O^{o_l} \circ E^{e_{l+1}})(n) = 1$ for positive integers $e_1, \dots, e_{l+1}$ and $o_1, \dots, o_l$. Let $\sigma_e = e_1+ \dots+ e_{l+1}$ and $\sigma_o = o_1+ \dots+ o_{l}$. Then, for $l>1$,
\[
\frac{n+1}{2^{\sigma_e}+1} < \left(\frac{2}{3}\right)^{\sigma_o}.
\]
\end{Corollary}
\begin{Proof}
The result follows from Theorem 1 and Lemma 4.
\end{Proof}
\vspace{.25 in}

\section {A generalization of the main results}

Now let us consider functions of the form
\[   
g(n) = 
\begin{cases}
Q(n) &\quad\text{if $P(n)$}\\
R(n) &\quad\text{if $\neg P(n)$}\\
\end{cases}
\]
where $n$ is a positive rational numbers, $P(n)$ is a logical predicate, functions $Q$ and $R$ map rational numbers to rational numbers and $[Q^j,R^k]$ is constant for any positive integers $j$ and $k$. Theorem 1 can be generalized to give us insight into the behavior of $g(n)$. Every repeated composition of $g(n)$ can be written as an equivalent composition of $Q$ and $R$ functions. It is trivial that there exists some $C$ such that
\[
(Q^{q_1} \circ R^{r_1} \circ \dotsb \circ Q^{q_n} \circ R^{r_n})(n) = (Q^{q_1 + \dotsb + q_n} \circ R^{r_1 + \dotsb + r_n})(n)+C
\]
for natural numbers $q_1,\dots,q_n$ and $r_1,\dots,r_n$. The nugget of insight here though is that by repeatedly applying the identity $[Q^j,R^k] \equiv Q^j \circ R^k - R^k \circ Q^j$ it is possible to iteratively build $(Q^{q_1 + \dotsb + q_n} \circ R^{r_1 + \dotsb + r_n})(n)$ from $(Q^{q_1} \circ R^{r_1} \circ \dotsb \circ Q^{q_n} \circ R^{r_n})(n)$ and acquire an expression for $C$ in the process. In some cases, as we have seen, the expression for $C$ derived in this way can be very manageable.

\section*{Statements and declarations}
We declare no conflicts of interest.

\bibliographystyle{abbrv}
\bibliography{root}

\begin{thebibliography}{1}

\bibitem{Barina}
D.~Barina.
\newblock Convergence verification of the collatz problem.
\newblock {\em The Journal of Supercomputing}, 77(3):2681--2688, 2021.

\bibitem{Feng}
G.~Feng.
\newblock Stability analysis of piecewise discrete-time linear systems.
\newblock {\em IEEE Transactions on Automatic Control}, 47(7):1108--1112, 2002.

\bibitem{Iserles}
A.~Iserles.
\newblock {\em A first course in the numerical analysis of differential
  equations}.
\newblock Number~44. Cambridge university press, 2009.

\bibitem{Lin}
Z.~Lin.
\newblock On asymptotic stabilizability of discrete-time linear systems with
  delayed input.
\newblock In {\em 2007 IEEE International Conference on Control and
  Automation}, pages 432--437. IEEE, 2007.

\bibitem{Mimuro}
T.~Mimuro et~al.
\newblock On certain simple cycles of the collatz conjecture.
\newblock {\em SUT. J. Math}, 37(2):79--89, 2001.

\end{thebibliography}
\end{document}